\documentclass[11pt, reqno]{amsart}
\usepackage{amsmath, amsthm, amscd, amsfonts, amssymb, graphicx, color}
\usepackage[bookmarksnumbered, colorlinks, plainpages]{hyperref}
\hypersetup{colorlinks=true,linkcolor=red, anchorcolor=green, citecolor=cyan, urlcolor=red, filecolor=magenta, pdftoolbar=true}
\newtheorem{theorem}{Theorem}[section]
\newtheorem{lemma}[theorem]{Lemma}
\newtheorem{proposition}[theorem]{Proposition}
\newtheorem{corollary}[theorem]{Corollary}
\theoremstyle{definition}
\newtheorem{definition}[theorem]{Definition}

\theoremstyle{remark}
\newtheorem{remark}[theorem]{Remark}
\numberwithin{equation}{section}

\title[Semi-Fredholm theory in $C^{*}-$algebras]{Semi-Fredholm theory in $C^{*}-$algebras}

\author[S. Ivkovi\'c]{Stefan Ivkovi\'c}
\address{Mathematical Institute of the Serbian Academy of Sciences and Arts,
	p.p. 367, Kneza Mihaila 36, 11000 Beograd, Serbia}
\email{stefan.iv10@outlook.com}



\date{\today}
\begin{document}

\maketitle

\begin{abstract}
	Ke\v{c}ki\'{c} and Lazovi\'{c} introduced an axiomatic approach to Fredholm theory by considering Fredholm type elements in a unital C*-algebra as a generalization of C*-Fredholm operators on the standard Hilbert C*-module introduced by Mishchenko and Fomenko and of Fredholm operators on a properly infinite von Neumann algebra introduced by Breuer. In this paper we establish semi-Fredholm theory in unital C*-algebras as a continuation of the approach by Ke\v{c}ki\'{c} and Lazovi\'{c}. We introduce the notion of semi-Fredholm type elements and semi-Weyl type elements. We prove that the difference of the set of semi-Fredholm elements and the set of semi-Weyl elements is open in the norm topology, that the set of semi-Weyl elements is invariant under perturbations by finite type elements, and several other results generalizing their classical counterparts. Also, we illustrate applications of our results to the special case of properly infinite von Neumann algebras and we obtain a generalization of the  punctured neighbourhood theorem in this setting. 
\end{abstract}

\baselineskip17pt

\section{Introduction}

The Fredholm and semi-Fredholm theory on Hilbert and Banach spaces started by studying the integral equations introduced in the pioneering work by Fredholm in 1903 in \cite{F}. After that, the abstract theory of Fredholm and semi-Fredholm operators on Hilbert and Banach spaces was further developed in numerous papers and books such as \cite{AP1}, \cite{AP2} and \cite{ZZRD}. 
In addition to classical semi-Fredholm theory on Hilbert and Banach spaces, several generalizations of this theory have been considered. Breuer for example started the development of Fredholm theory in von-Neumann algebras as a generalization of the classical Fredholm theory for operators on Hilbert spaces. In \cite{BR} and \cite{BR2} he introduced the notion of a Fredholm operator in a von Neumann algebra and established its main properties.  On the other hand, Fredholm theory on Hilbert $C^{*}$-modules as another generalization of the classical Fredholm theory on Hilbert spaces was started by Mishchenko and Fomenko. In \cite{MF} they introduced the notion of a Fredholm operator on the standard Hilbert  $C^{*}$-module and proved a generalization in this setting of some of the main results from the classical Fredholm theory. In \cite{IS1}, \cite{IS3}, \cite{IS5} and \cite{IS6} we went further in this direction and defined semi-Fredholm and semi-Weyl operators on Hilbert $C^{*}$-modules. We investigated and proved several properties of these new semi- Fredholm operators on Hilbert $C^{*}$-modules as a generalization of the results from the classical semi-Fredholm theory on Hilbert and Banach spaces. 
The interest for considering these generalizations comes from the theory of pseudo differential operators acting on manifolds. The classical theory can be applied in the case of compact manifolds, but not in the case of non-compact ones. Even operators on Euclidian spaces are hard to study, for example Laplacian is not Fredholm. Kernels and cokernels of many operators are infinite dimensional Banach spaces, however, they may also at the same time be finitely generated Hilbert modules over some appropriate $C^{*}$-algebra. Similarly, orthogonal projections onto kernels and cokernels of many bounded linear operators on Hilbert spaces are not finite rank projections in the classical sense, but they are still finite projections in some appropriate von Neumann algebra. Therefore, many operators that are not semi-Fredholm in the classical sense may become semi-Fredholm in a more general sense if we consider them as operators on Hilbert $C^{*}$-modules or as elements of an appropriate von Neumann algebra. Hence, by studying these generalized semi-Fredholm operators, we get a proper extension of the classical semi-Fredholm theory to new classes of operators.

Now, Ke\v{c}ki\'{c} and Lazovi\'{c} in \cite{KL} established an axiomatic approach to Fredholm theory. They introduced the notion of a finite type element in a unital $C^{*}$-algebra which generalizes the notion of the compact operator on the standard Hilbert $C^{*}$-module and the notion of a finite operator in a properly infinite von Neumann algebra. They also introduced the notion of a Fredholm type element with respect to the ideal of these finite type elements. This notion is at a same time a generalization of the classical Fredholm operator on a Hilbert space, Fredholm $C^{*}$-operator on the standard Hilbert $C^{*}$-module defined by Mishchenko and Fomenko and the Fredholm operator on a properly infinite von Neumann algebra defined by Breuer. The index of this Fredholm type element takes values in the K-group. They obtained that the set of Fredholm type elements in a unital $C^{*}$-algebra is open in the norm topology and they proved a generalization of the Atkinson theorem. Moreover, they proved the multiplicativity of the index in the K-group and that the index is invariant under perturbation of Fredholm type elements by finite type elements.   

The aim of this paper is to establish semi-Fredholm theory in unital $C^{*}$-algebras as a continuation of the approach by Ke\v{c}ki\'{c} and Lazovi\'{c} on Fredholm theory in unital $C^{*}$-algebras. We introduce the notion of a semi-Fredholm type element and a semi-Weyl type element with respect to the ideal of finite type elements and obtain a generalization in this setting of several results from the classical semi-Fredholm and semi- Weyl theory of operators on Hilbert spaces. The motivation for this research is not only developing an abstract, axiomatic semi-Fredholm theory in unital $C^{*}$-algebras, but also to obtain an extension of Breuer`s Fredholm theory to semi-Fredholm and semi-Weyl theory in properly infinite von Neumann algebras by applying our results to this special case. In the third section of the paper, we present the results in abstract semi-Fredholm theory and semi-Weyl theory in unital $C^{*}$-algebras, whereas in the fourth section of the paper we give the application of these results to the concrete case of properly infinite von Neumann algebras. 

The paper is organized as follows. \\
In the third section, after introducing some technical lemmas and corollaries, we show that an element in a unital  $C^{*}$-algebra $\mathcal{A}$ is of Fredholm type if and only if its compressions with respect to finite type projections are of Fredholm type, and in this case the indexes in the K-group are the same. We define semi-Fredholm type elements and provide several characterizations of these elements including a generalization in this setting of the classical Fredholm alternative. Also, we define semi-Weyl type elements and prove that the set of these elements and the difference of the set of semi-Fredholm type elements and the set of semi-Weyl type elements are both open in the norm topology of $\mathcal{A}$. This result has several consequences in connection with continuous maps from the unit interval into the set of semi-Fredholm type elements and with topological properties of the set of semi-Fredholm type elements, as illustrated in corollaries. Next, we show the set of semi-Weyl type elements is invariant under perturbations by finite type elements. This allows us to give an algebraic description of semi-Weyl type elements in terms of left and right invertible elements in $\mathcal{A},$ from which we deduce that the set of semi-Weyl type elements is closed under the multiplication. Finally, we consider generalized Weyl type elements in $\mathcal{A}$ as a generalization of generalized Weyl operators on Hilbert spaces defined by \DJ{}or\dj{}evi\'{c} in \cite{DDj2}, and we extend in this setting \DJ{}or\dj{}evi\'{c}s main result in \cite{DDj2} by proving that a product of two generalized Weyl type elements in $\mathcal{A}$ is also a generalized Weyl type element provided that this product admits Moore-Penrose inverse in $\mathcal{A}.$ In the special case when $\mathcal{A}$ is a von Neumann algebra, this result gives us that if two operators belonging to $\mathcal{A}$ satisfy that their kernels are Murray von Neumann equivalent to the orthogonal complement of their images and if their composition has closed image, then the kernel of their composition is Murray von Neumann equivalent to the orthogonal complement of the image of their composition.\\
In the fourth section of the paper, we establish semi-Fredholm and semi-Weyl theory in properly infinite von Neumann algebras as a continuation Breuer`s Fredholm theory in von Neumann algebras. Among other results we show that an operator belonging to $\mathcal{A}$ is upper (respectively lower) semi-Weyl in this new generalized sense if and only if it can be written as a sum of a bounded below (respectively surjective) operator and a finite operator. Also, we show that the set of those Fredholm operators which are both upper semi-Weyl and lower semi-Weyl is exactly the set of Weyl operators in $\mathcal{A},$ and we provide a generalization of the punctured neighbourhood theorem in the setting of Breuer`s Fredholm theory in von Neumann algebras.

\section{Preliminaries }

In this section we recall the following definitions.

\begin{definition} \label{d 09}
	\cite[Definition 1.1]{KL} Let $\mathcal{A} $ be an unital $C^{*}$-algebra, and $\mathcal{F} \subseteq \mathcal{A} $ be a subalgebra which satisfies the following conditions:\\
	(i) $\mathcal{F} $ is a selfadjoint ideal in $\mathcal{A} ,$ i.e. for all $a \in \mathcal{A}, b \in \mathcal{F}  $ there holds $ab,ba \in \mathcal{F} ,$ and $a \in \mathcal{F} $ implies $a^{*} \in \mathcal{F} ;$\\
	(ii) There is an approximate unit $p_{\alpha}$ in the norm topology for $\mathcal{F}$ consisting of projections;\\
	(iii) If $p,q \in \mathcal{F} $ are projections, then there exists $v \in \mathcal{A} ,$ such that $vv^{*}=q $ and $v^{*}v \perp p, $ i.e. $v^{*}v + p $ is a projection as well;\\
	The elements of such ideal we shall call \textit{finite type elements}. In further, we shall denote it by $\mathcal{F} .$
\end{definition}

\begin{definition} \label{d 10} \cite[Definition 1.2]{KL}
	Let $\mathcal{A} $ be an unital $C^{*}-$ideal, and let $\mathcal{F} \subseteq \mathcal{A}$ be an algebra of finite type elements.In the set $ \text{ Proj}(\mathcal{F}) $ we define the equivalence relation:  
	$$p \sim q \Leftrightarrow \exists v \in \mathcal{A} \text{ } vv^{*}=p,\text{ } v^{*}v=p, $$
	i.e. Murray - von Neumann equivalence. The set $ S(\mathcal{F})=\text{Proj}(\mathcal{F})\text{ }/\sim $ is a commutative semigroup with respect to addition, and the set,$K(\mathcal{F})=G(S(\mathcal{F})),$ where $G$ denotes the Grothendic functor, is a commutative group.
\end{definition}

\begin{definition} \label{d 11} \cite[Definition 2.1]{KL}
	Let $a \in \mathcal{A} $ and $p,q $ be projections in $\mathcal{A} .$ Set $a^{\prime}=(1-q)a(1-p) .$ We say that $a$ is invertible up to pair $(p, q)$ if there is some $b \in \mathcal{A} $ with $b=(1-p)b(1-q) $ (and immediately $bp=0, $ $pb=0,$ $b=(1-p)b=b(1-q )$) such that 
	$$a^{\prime}b=1-q, \text{  } ba^{\prime}=1-p. $$
	We refer to such $b$ as almost inverse of $a,$ or $(p,q)$-inverse of $a.$
\end{definition}

\begin{lemma} \label{l 01}
	Let $a \in \mathcal{A} $ and $p,q $ be projections in $\mathcal{A} .$ Then $a$ is invertible up to pair $(p,q) $ if and only if $a^{*}$ is invertible up to pair $(q,p) .$ 
\end{lemma}

\begin{proof}
	We just observe that $a^{\prime} b= 1-q $ and $ ba^{\prime}$ if and only if $b^{*}a^{\prime *} =1-q $ and $a^{\prime *}b^{*}=1-p   ,$ where $a^{\prime} = (1-q)a(1-p) .$ Since $a^{\prime *} = (1-p)a^{*}(1-q) ,$ the statement follows. 
\end{proof}

\begin{definition} \label{d 12} \cite[Definition 2.2]{KL}
	Let $\mathcal{F}$ be the ideal of finite type elements. We say that $a \in \mathcal{A} $ is of Fredholm type (or abstract Fredholm element) with respect to the ideal $\mathcal{F}$ if there are projections $p,q \in \mathcal{F} $ such that $a$ is invertible up to $(p,q).$ The index of the element $a$ (or abstract index) is the element of the group $K(\mathcal{F}) $ defined by
	$$\text{ind}(a)=([p],[q]) \in K(\mathcal{F}), $$
	or less formally
	$$\text{ind}(a)=[p]-[q]. $$
\end{definition}

\section{Main results}

Throughout this paper, $\mathcal{A}$ always stands for a unital $C^{*}$-algebra and $\mathcal{F}$ denotes ideal of finite type elements.

The following lemma is a technical result which is needed in several applications later in the paper. 

\begin{lemma} \label{l 02}
	Let $a \in \mathcal{A} $ and $p,q,p^{\prime},q^{\prime} $ be projections in $\mathcal{A} .$ Suppose that $p,q,p^{\prime} \in \mathcal{F} .$ If $a$ is invertible up to pair $(p,q) $ and also invertible up to pair $(p^{\prime},q^{\prime}) ,$ then $q^{\prime} \in \mathcal{F} .$ 
	Similarly, if instead of $p,q,p^{\prime} $ we have that $p,q,q^{\prime} \in \mathcal{F} ,$ then we must have that $p^{\prime} \in \mathcal{F}$ as well.
\end{lemma}

\begin{proof}
	By \cite[Proposition 2.8]{KL} we may without loss of generality assume that $qa(1-p)=0 ,$ hence we can apply \cite[Proposition 2.11]{KL} with respect to the pair $(p,q) .$ Let $p_{\alpha} $  and $q_{\alpha}$ be the approximate units from \cite[Proposition 2.11]{KL} with respect to the pair $(p,q) .$ By \cite[Lemma 2.4]{KL}, there is some $\delta >0 $ such that if $\parallel p^{\prime \prime} - p^{\prime} \parallel < \delta  $ for some projection $ p^{\prime \prime} $ in $\mathcal{A},$ then $a$ is also invertible up to pair $(p^{\prime \prime}, q^{\prime} ).$ Now, since $p_{\alpha} $ is an approximate unit for  $\mathcal{F} ,$ by \cite[Lemma 2.9]{KL} there is some $\alpha_{0} $ and some projection $p^{\prime \prime} $ in $\mathcal{A} $ such that $p^{\prime} \sim  p^{\prime \prime} \leq p_{\alpha} $ for all $\alpha \geq \alpha_{0} $ and $\parallel p^{\prime \prime} - p^{\prime} \parallel < \delta .$ Since $a$ is then invertible up to pair $( p^{\prime \prime}, q^{\prime} ),$ by  \cite[Lemma 2.9]{KL} there is a projection $q^{\prime \prime} $ in $\mathcal{A} $ such that $ q^{\prime \prime} \sim q^{\prime},$ $a$ is invertible up to pair $(p^{\prime \prime}, q^{\prime \prime}) $ and $q^{\prime \prime} a (1-p^{\prime \prime})=0 .$ Hence we can repeat the proof of \cite[Proposition 2.11]{KL} in order to construct an a net of projections $\lbrace q_{\alpha}^{\prime \prime} \rbrace_{\alpha \geq \alpha_{0}} $ such that the pairs $(p_{\alpha} , q_{\alpha}^{\prime \prime})_{\alpha \geq \alpha_{0}} $ satisfy the conditions of \cite[Proposition 2.11]{KL} with respect to the pair $(p^{\prime \prime}, q^{\prime \prime}) .$ 
	Now, by \cite[Proposition 2.8]{KL}, for each $\alpha \geq \alpha_{0} $ there exist projections $\tilde{q}_{\alpha} $ and $\tilde{q}_{\alpha}^{\prime \prime} $ such that $\tilde{q}_{\alpha} \sim q_{\alpha},$ $\tilde{q}_{\alpha}^{\prime \prime} \sim {q}_{\alpha}^{\prime \prime} ,$ $\tilde{q}_{\alpha} a (1-p_{\alpha})=\tilde{q}_{\alpha}^{\prime \prime} a (1-p_{\alpha})=0 ,$ and $a$ is invertible both up to pair $(p_{\alpha} , \tilde{q}_{\alpha} ) $ and pair $(p_{\alpha} , \tilde{q}_{\alpha}^{\prime \prime}  ) .$ Hence,
	$$a(1-p_{\alpha})=(1-\tilde{q}_{\alpha})a(1-p_{\alpha})=(1-\tilde{q}_{\alpha}^{\prime \prime}) a (1-p_{\alpha}) $$ 
	for all $\alpha \geq \alpha_{0} .$ Moreover, $a$ is invertible both up to pair $(p_{\alpha},\tilde{q}_{\alpha}) $ and $(p_{\alpha} ,\tilde{q}_{\alpha}^{\prime \prime}).$ Let $b_{\alpha} $ be $(p_{\alpha} , \tilde{q}_{\alpha})$ inverse of $a$ for $\alpha \geq \alpha_{0} .$ Then we have 
	$$ 0 = \tilde{q}_{\alpha}^{\prime \prime} a (1-p_{\alpha}) b_{\alpha} = \tilde{q}_{\alpha}^{\prime \prime} (1-\tilde{q}_{\alpha}) a (1-p_{\alpha}) b_{\alpha} = \tilde{q}_{\alpha}^{\prime \prime} (1-\tilde{q}_{\alpha}).$$ 
	It follows that $\tilde{q}_{\alpha}^{\prime \prime} \leq \tilde{q}_{\alpha} $ for all $\alpha \geq \alpha_{0} .$ Similarly we can show that $ \tilde{q}_{\alpha} \leq \tilde{q}_{\alpha}^{\prime \prime},$ hence $\tilde{q}_{\alpha} = \tilde{q}_{\alpha}^{\prime \prime} .$ Thus $ {q}_{\alpha}^{\prime \prime} \sim \tilde{q}_{\alpha}^{\prime \prime} \sim {q}_{\alpha} \in \mathcal{F}$ for all $\alpha \geq \alpha_{0} .$ However, we have by the construction of $q_{\alpha}^{\prime \prime} $ that $ q_{\alpha}^{\prime \prime} - q^{\prime \prime} \sim p_{\alpha} - p$ which is an element of $\mathcal{F}$. Hence we must have $q^{\prime \prime} \in \mathcal{F} .$ Since $q^{\prime \prime} \sim q^{\prime} ,$ we get that $q^{\prime} \in \mathcal{F} .$\\
	The second statement can be proved by passing to the adjoints and applying Lemma \ref{l 01}.
	
\end{proof}

The following corollary is a generalization of Fredholm alternative. 

\begin{corollary} \label{l 08}
	Let $a \in G(\mathcal{A}) $ and suppose that $K(\mathcal{F}) $ satisfies the cancellation property i.e. for any pair of projections $p,q $ in $\mathcal{F} $ we have that $p \sim q $ whenever $ [p]=[q].$ Then for every $f \in \mathcal{F} $ we have that $ a+f$ is left invertible in $\mathcal{A} $ if and only if $ a+f $ is right invertible in $\mathcal{A} .$ 
\end{corollary}

\begin{proof}
	Since $a$ is invertible by assumption, we have that $a$ is of Fredholm type and  $index \text{ }a=0$ in $ K(\mathcal{F}).$ Next, by part $a)$ in \cite[Proposition 2.14]{KL} we have index $index \text{ }(a+f)=index \text{ }a=0.$ If $a+f$ is left invertible, by \cite[Lemma 2.5]{KL} there exists a projection $r \in \mathcal{A} $ such that $(1-r)(a+f)=a+f $ and such that $a+f$ is invertible up to pair $(0,r) .$ By Lemma \ref{l 02} it follows that $r \in \mathcal{F} .$ Hence $0 = index \text{ } (a+f)=-[r] .$ Since $K(\mathcal{F}) $ satisfies the cancellation property by assumption, we get that $r=0 ,$ so $a+f$ is invertible. By passing to the adjoints and applying Lemma \ref{l 01} we can deduce that $a+f$ is invertible in $\mathcal{A} $ if $a+f$ is right invertible in $\mathcal{A} .$
\end{proof}

From the proof of Lemma \ref{l 02} we can also deduce the following corollaries. 

\begin{corollary} \label{cor 3.3}
	Let $ a \in \mathcal{A} $ and $p,q, p^{\prime}, q^{\prime}$ be projections in $\mathcal{A}.$  If $a$ is invertible both up to pair $(p,q)$ and up to pair $(p , q^{\prime}),$ then $q \sim q^{\prime} .$ If $a$ is invertible both up to pair $(p,q)$ and $(p^{\prime},q),$ then $p \sim p^{\prime}.$
\end{corollary}

Now we provide a special characterization of Fredholm type elements which we will use later when we consider the compressions of elements in $\mathcal{A}$ with respect to finite type projections. 


\begin{lemma} \label{l 03}
	Let $ a \in \mathcal{A} .$ Then $a$ is a Fredholm type element with respect to the ideal $\mathcal{F}$ if and only if for every projection $ \tilde{p} \in \mathcal{F}$ there exist some projections $p^{\prime}, q^{\prime} \in \mathcal{F} $ such that $\tilde{p} \leq p^{\prime} ,$ $\tilde{p} \leq q^{\prime} $ and $a$ is invertible up to pair $(p^{\prime}, q^{\prime}) .$
\end{lemma}

\begin{proof}
	Since $\tilde{p} \in \mathcal{F} ,$ we have that $1-\tilde{p} $ is a Fredholm type element. If $a$ is also a Fredholm type element, then $(1-\tilde{p})  a (1-\tilde{p}) $ is a Fredholm type element as well. Hence there exist some projections $p,q \in \mathcal{F} $ such that $(1-\tilde{p})  a (1-\tilde{p})$ is invertible up to pair $(p,q).$ In exactly the same way as in the proof of the part $b)$ in \cite[Proposition 2.15]{KL}, we can deduce that there exists a projection $r \in \mathcal{A} $ such that $(1-r) (1-\tilde{p}) (1-p) = (1-\tilde{p}) (1-p) $ and such that $1 - \tilde{p} $ is invertible up to pair $(p, r.)$ Now, since $1-\tilde{p} $ is a Fredholm type element, by Lemma \ref{l 02} it follows that $r \in \mathcal{F}.$ Let $c$ be $(p, r)$ inverse of $1-\tilde{p} .$ Then we have 
	$$\tilde{p} (1-r) = \tilde{p} (1-r) (1-\tilde{p}) (1-p) c =  \tilde{p}  (1-\tilde{p}) (1-p) c = 0   .$$ 
	Hence $ \tilde{p} \leq r .$ Next we show that $(1-\tilde{p}) a$ is invertible up to pair $(r, q).$ Indeed if $b$ is $(p, q)$ inverse of $(1-\tilde{p}) a (1-\tilde{p}),$ then we have 
	$$(1-\tilde{p}) (1-p) b (1-q) (1-\tilde{p}) a (1-r) =(1-\tilde{p}) (1-p) b (1-q) (1-\tilde{p}) a (1-r) (1-\tilde{p}) (1-p) c $$
	$$=(1-\tilde{p}) (1-p) b (1-q) (1-\tilde{p}) a  (1-\tilde{p}) (1-p) c =(1-\tilde{p}) (1-p) (1-p) c= (1-\tilde{p}) (1-p) c $$
	$$=(1-r)(1-\tilde{p}) (1-p) c = 1-r $$ 
	Moreover, 
	$$(1-q) (1-\tilde{p}) a (1-r) (1-\tilde{p}) (1-p) b = (1-q) (1-\tilde{p}) a (1-\tilde{p}) (1-p) b = 1-q .$$ 
	Hence $(1-\tilde{p})(1-p)b$ is $(r,q)$ inverse of $(1-\tilde{p})a.$
	By \cite[Proposition 2.8]{KL} there exists some projection $q^{\prime} \sim q $ such that $(1-\tilde{p})a $ is invertible up to pair $(r, q^{\prime}) $ and $ q^{\prime}(1-\tilde{p}) a (1-r) =0 .$ Let $d$ be $ (r, q^{\prime}) $ inverse $(1- \tilde{p})a .$ Then we get $$\tilde{p} (1-q^{\prime}) = \tilde{p} (1-q^{\prime}) (1-\tilde{p}) a (1-r) d=\tilde{p}  (1-\tilde{p}) a (1-r)d = 0   .$$ 
	Hence, $ \tilde{p} \leq q^{\prime}.$ Furthermore, $(1-q^{\prime}) (1-\tilde{p}) = (1-q^{\prime})  ,$ so $a$ is invertible up to pair $(r, q^{\prime}) ,$ which proves the implication in one direction. The implication in the opposite direction follows from the definition of Fredholm type element. 
\end{proof} 

From the proof of Lemma \ref{l 03} we can deduce the following corollary. 

\begin{corollary} \label{c 05}
	Let $a,b \in \mathcal{A} $ and $p, q$ be projections in $\mathcal{A}$. Suppose that $ba$ is invertible up to pair $(p, q).$ Then there exists a projection $r \in \mathcal{A} $ such that $a$ and $b$ are invertible up to pairs $(p, r)$ and $(r, q),$ respectively. Moreover, we have that $(1-r) a (1-p) = a (1-p) .$ 
\end{corollary}

\begin{lemma}
	Let $p \in \mathcal{F} $ be a projection. Then the couple 
	$$(  (1-p)\mathcal{A} (1-p) ,  (1-p)\mathcal{F} (1-p) ) $$ 
	satisfies the conditions of Definition \ref{d 09}.
\end{lemma}

\begin{proof}
	It is obvious that $ (1-p) \mathcal{F} (1-p)$ is a self-adjoint ideal in the corner $C^{*}$-algebra $(1-p) \mathcal{A} (1-p) $ since $\mathcal{F} $ is a self-adjoint ideal in $\mathcal{A} .$ Next, 
	by the straightforward calculation one can show that
	$\lbrace (1-p) p_{\alpha} (1-p) \rbrace$ is an a approximate unit in\\
	$(1-p)\mathcal{F} (1-p) .$ 
	
	Futher, let $\tilde{p} , \tilde{q} $ be projections in $ (1-p)\mathcal{F} (1-p)  .$ Since $ (1-p) \tilde{p} (1-p)= \tilde{p}  ,$ it follows that $\tilde{p} \perp p .$ By the property $(iii)$ in {d 09}, since $p+ \tilde{p} $ and $\tilde{q} $ are projections in $ \mathcal{F},$ there exists some $ v \in \mathcal{A}$ such that $vv^{*}= \tilde{q} $ and $v^{*}v $ is a projection satisfying that $ v^{*}v \perp p+\tilde{p}.$ Set $q^{\prime}= v^{*}v .$ Then 
	$( \tilde{q} v q^{\prime} ) (\tilde{q} v q^{\prime} )^{*}=\tilde{q} $ and 
	$( \tilde{q} v q^{\prime} )^{*} (\tilde{q} v q^{\prime} )=q^{\prime}  .$ However, since $q^{\prime} \perp p+\tilde{p} ,$ we have that $\tilde{q}, q^{\prime} \in (1-p)\mathcal{A} (1-p) .$ Moreover, in particular $q^{\prime} \perp p ,$ hence the couple $(  (1-p)\mathcal{A} (1-p) ,  (1-p)\mathcal{F} (1-p) ) $ satisfies the condition $(iii)$ in Definition \ref{d 09}.
\end{proof}

Now we can relate Fredholmness of elements in $\mathcal{A}$ to Fredholmness of the compressions of these elements with respect to finite type projections, thus generalizing in this setting \cite[Lemma 2.10.1]{ZZRD} originally given in \cite{PO}.

\begin{corollary}
	Let $a \in \mathcal{A} $  and $p $ be a projection in $\mathcal{F} .$ Then $a$ is a Fredholm type element in $ \mathcal{A}$ with respect to the ideal $\mathcal{F}$ if and only if $ (1-p) a (1-p) $ is a Fredholm type element in $(1-p) \mathcal{A} (1-p) $  with respect to the ideal $(1-p) \mathcal{F} (1-p),$ and in this case $index \text{ } a=index \text{ } (1-p) a (1-p).$ 
\end{corollary}

\begin{proof}
	If $a$ is a Fredholm type element in $\mathcal{A} ,$ then by Lemma \ref{l 03} there exist projections $p^{\prime}, q^{\prime} \in \mathcal{F} $ such that $p \leq p^{\prime}, p \leq q^{\prime} $ and $a$ is invertible up to $(p^{\prime},q^{\prime}) .$ It is straightforward to see that $(1-p) a (1-p)$ is invertible in $(1-p) \mathcal{A} (1-p) $ up to $(p^{\prime}-p,q^{\prime}-p) .$ Indeed, recall that $1-p$ is the unit in $(1-p) \mathcal{A} (1-p) $ and if $b$ is $(p^{\prime}, q^{\prime}) $ inverse of $a,$ then $(1-p^{\prime}) b (1-q^{\prime}) $ is also $(p^{\prime} , q^{\prime}) $ inverse of $a.$
	Since $p^{\prime}-p,q^{\prime}-p \in (1-p) \mathcal{F} (1-p) ,$ it follows that $(1-p) a (1-p)$ is Fredholm type element in $(1-p) \mathcal{A} (1-p) $ and $index \text{ } (1-p) a (1-p) = [p^{\prime}]-[q^{\prime}] = index \text{ } a .$ Conversely, if $(1-p) a (1-p)$ is Fredholm type element in $(1-p) \mathcal{A} (1-p) ,$ then there exist projections $\tilde{p}, \tilde{q} \in  (1-p) \mathcal{F} (1-p)$ such that $(1-p) a (1-p)$ is invertible in $(1-p) \mathcal{A} (1-p) $ up to $ (\tilde{p}, \tilde{q}).$ It follows that $a$ is invertible in $\mathcal{A} $ up to $(p+\tilde{p},p+\tilde{q}) $ and $index \text{ }a=[\tilde{p}]-[\tilde{q}] = index \text{ } (1-p) a (1-p).$
\end{proof}

We define now semi-Fredholm elements in $\mathcal{A}.$

\begin{definition} \label{d 04}
	Let $a \in \mathcal{A} .$ We say that $a$ is an upper semi-Fredholm element with respect to the ideal $\mathcal{F}$ if $a$ is invertible up to pair of projections $(p, q)$ where $p \in \mathcal{F} .$ Similarly, we say that $a$ is a lower semi-Fredholm element with respect to the ideal $\mathcal{F}$, however in this case we assume that $q \in \mathcal{F} $ (and not $p$).
\end{definition}

The next two lemmas give some characterizations of semi-Fredholm elements in $\mathcal{A}.$ 

\begin{lemma} \label{l 06}
	Let $a \in \mathcal{A} .$ Then the following holds.\\
	1) If $a$ is an upper semi-Fredholm element and $\tilde{p} $ is a projection in $\mathcal{A} $ such that $a \tilde{p}=0, ,$ then $\tilde{p} \in \mathcal{F} .$ In particular, if $p, q$ are projections in $\mathcal{A} $ such that $a$ is invertible up to $(p, q)$ and $(1-q) a (1-p)=a,$ then $a$ is an upper semi-Fredholm element if and only if $p \in \mathcal{F} .$\\
	2) If $a$ is a lower semi-Fredholm element and $\tilde{q} $ is a projection in $\mathcal{A} $ such that $\tilde{q} a = 0 ,$ then $ \tilde{q} \in \mathcal{F}.$ In particular, if $p, q$ are projections in $\mathcal{A} $ such that $a$ is invertible up to $(p, q)$ and $(1-q) a (1-p)=a,$ then $a$ is a lower semi-Fredholm element if and only if $q\in \mathcal{F} .$
\end{lemma}

\begin{proof}
	We prove first $1).$ Suppose that $a$ is invertible up to pair $(p, q)$ where $p \in \mathcal{F} .$ By \cite[ Proposition 2.8]{KL} we may without loss of generality assume that $(1-q)ap=0 .$ Then we have 
	$$ (1-p)\tilde{p}=b(1-q)a(1-p)\tilde{p} = b (1-q)a\tilde{p}=0$$ where $b$ is $(p, q)$ 
	inverse of $a.$ Hence we get that $\tilde{p}=p\tilde{p} $ which is an element of $ \mathcal{F} $ since $p \in \mathcal{F} $  and $ \mathcal{F} $ is an ideal in $ \mathcal{A}.$ By passing to the adjoints and applying Lemma \ref{l 01} we can prove $2).$
\end{proof}

\begin{lemma} \label{l 07}
	Let $a \in \mathcal{F} .$ Then $a$ is an upper semi-Fredholm element if and only if $a$ is left invertible up to some projection $p \in \mathcal{F} .$ Similarly, $a$ is a lower semi-Fredholm element if and only if $a$ is right invertible up to some projection $q \in \mathcal{F} .$ 
\end{lemma}

\begin{proof}
	Suppose that $a$ is left invertible up to some projection $p \in \mathcal{F} .$ Then, by \cite[Lemma 2.5]{KL} there exists a projection $ r \in \mathcal{A}$ such that $(1-r)a(1-p)=a(1-p) $ and such that $a$ is invertible up to pair $(p, r).$ Hence $a$ is an upper semi-Fredholm element then. Conversely, if $a$ is upper semi-Fredholm element, then $a$ is invertible up to some pair $(p, q)$ where $p \in \mathcal{F} .$ Hence $a$ is left invertible up to $p.$\\
	By passing to the adjoints and applying Lemma \ref{l 01} we can deduce the second statement. 
\end{proof}

\begin{definition} \cite{BR} \cite{BR2} \label{r10d 1.1}
	Let $p, q$ be projections in $\mathcal{A}.$ We will denote $p \preceq q $ if there exists some projection $p^{\prime} $ such that $p^{\prime} \leq q $ and $ p \sim p^{\prime} .$ 
\end{definition}

We recall now the following auxiliary technical lemma.

\begin{lemma}    \label{r10l 1.2}
	1) If $p_{1} \sim q_{1},$ $p_{2} \sim q_{2}$ and $p_{1}p_{2}=q_{1}q_{2}=0 ,$ then $(p_{1}+p_{2}) \sim (q_{1}+q_{2}) .$\\
	2) If $p \preceq q $ and $q \sim q^{\prime} ,$ then $p \preceq q^{\prime} .$ 
\end{lemma} 

\begin{corollary} \label{r10c 1.5}
	Let $ p_{1}, p_{2}, q_{1}, q_{2} $ be projections in $\mathcal{A}$ such that $ p_{1} p_{2}=q_{1} q_{2}=0, \text{ }  p_{1} \sim q_{1}$ and $p_{2} \preceq q_{2} .$ Then $p_{1} + p_{2} \preceq q_{1} + q_{2} .$ 
\end{corollary} 

Now we introduce the notion of semi-Weyl type elements in $\mathcal{A}$ as a generalization of semi-Weyl operators on Hilbert spaces. 

\begin{definition} \label{r10d 1.3}
	Let $a \in \mathcal{A} .$ We say that $a$ is an upper semi-Weyl type element with respect to the ideal $\mathcal{F}$ if there exist projections $p, q$ in $\mathcal{A} $ such that $p \in \mathcal{F} ,$ $p \preceq q $ and $a$ is invertible up to pair $(p, q).$ Similarly we say that $a$ is a lower semi-Weyl type element with respect to the ideal $\mathcal{F}$, only in this case we assume that $q \in \mathcal{F} $ and $q \preceq p .$ Finally, we say that $a$ is a Weyl type element with respect to the ideal $\mathcal{F}$ if $a$ is invertible up to pair $(p, q)$ where $p, q $ are projections in $ \mathcal{F} $ and $p \sim q .$
\end{definition}

Set
\begin{flushleft}
	$$\mathcal{K}\Phi_{+} (\mathcal{A})= \lbrace a \in \mathcal{A} \mid  a \text{ is upper semi-Fredholm type element }  \rbrace ,$$
	$$\mathcal{K}\Phi_{-} (\mathcal{A})= \lbrace a \in \mathcal{A} \mid  a \text{  is lower semi-Fredholm type element }  \rbrace ,$$
	$$\mathcal{K}\Phi (\mathcal{A})= \lbrace a \in \mathcal{A} \mid  a \text{ is  Fredholm type element }  \rbrace ,$$
	$$\mathcal{K}\Phi_{+}^{-} (\mathcal{A})= \lbrace a \in \mathcal{A} \mid  a \text{  is upper semi-Weyl type element }  \rbrace ,$$
	$$\mathcal{K}\Phi_{-}^{+} (\mathcal{A})= \lbrace a \in \mathcal{A} \mid a \text{  is lower semi-Weyl type element }  \rbrace ,$$
	$$\mathcal{K}\Phi_{0} (\mathcal{A})= \lbrace a \in \mathcal{A} \mid  a \text{  is Weyl type element }  \rbrace .$$
\end{flushleft}
It is understood that we consider a fixed ideal $\mathcal{F}$ of finite type elements. \\
Notice that by definition we have $\mathcal{K}\Phi_{+}^{-} (\mathcal{A}) \subseteq \mathcal{K}\Phi_{+} (\mathcal{A}), \text{ } \mathcal{K}\Phi_{-}^{+} (\mathcal{A}) \subseteq \mathcal{K}\Phi_{-} (\mathcal{A}) $ and $\mathcal{K}\Phi_{0} (\mathcal{A}) \subseteq \mathcal{K}\Phi (\mathcal{A}) .$ 

\begin{remark}
	In exactly the same way as in the proof of \cite[Proposition 2.15]{KL} one can show that $a \in\mathcal{K}\Phi_{+} (\mathcal{A}) $ if and only if $a$ is left invertible modulo $\mathcal{F} $ and similarly, $a \in\mathcal{K}\Phi_{-} (\mathcal{A}) $ if and only if $a$ is right invertible modulo $\mathcal{F} .$ Hence, an analogue in this setting of \cite[Corollaries 2.4-2.10]{IS1} hold. From Lemma \ref{l 01} it follows that $a \in\mathcal{K}\Phi_{+} (\mathcal{A}) $ if and only if $a^{*} \in\mathcal{K}\Phi_{-} (\mathcal{A}) $ and similarly, $ a \in\mathcal{K}\Phi_{+}^{-} (\mathcal{A})$ if and only if $a^{*} \in\mathcal{K}\Phi_{-}^{+} (\mathcal{A}) .$ Moreover, by \cite[Lemma 2.4]{KL} the sets $\mathcal{K}\Phi_{+} (\mathcal{A}), \mathcal{K}\Phi_{-} (\mathcal{A}), \mathcal{K}\Phi_{+}^{-} (\mathcal{A}) , \mathcal{K}\Phi_{-}^{+} (\mathcal{A})  $ and $\mathcal{K}\Phi_{0} (\mathcal{A})  $ are open in the norm topology of $\mathcal{A} .$ 
\end{remark}

\begin{proposition} \label{r10p 1.4}
	The sets $\mathcal{K}\Phi_{+}(\mathcal{A}) \setminus \mathcal{K}\Phi_{+}^{-} (\mathcal{A})  ,$ $\mathcal{K}\Phi_{-}(\mathcal{A}) \setminus \mathcal{K}\Phi_{-}^{+} (\mathcal{A})  $ and $\mathcal{K}\Phi(\mathcal{A}) \setminus \mathcal{K}\Phi_{0} (\mathcal{A}) $ are open in the norm topology of $\mathcal{A} .$
\end{proposition} 

\begin{proof}
	Let $a \in \mathcal{K}\Phi_{+}(\mathcal{A}) \setminus \mathcal{K}\Phi_{+}^{-} (\mathcal{A}) $ and $p, q$ be the projections in $\mathcal{A} $ such that $p \in \mathcal{F} $ and such that $a$ is invertible up to pair $(p, q).$ By \cite[Proposition 2.8]{KL} we may without loss of generality assume that $qa(1-p)=0 .$ By \cite[Lemma 2.4]{KL}, there exists an $\epsilon > 0 $ such that if $\parallel a-b \parallel < \epsilon $ for some $b \in \mathcal{A} ,$ then $b $ is also invertible up to pair $(p, q).$ Now, by \cite[Proposition 2.8]{KL}, there exists some projection $q^{\prime} $ in $\mathcal{F} $ such that $q \sim q^{\prime} ,$ $b $ is invertible up to pair $(p, q^{\prime}) $ and $q^{\prime}b(1-p)=0 .$ Suppose that there exist some projections $p^{\prime}, q^{\prime \prime}  \in \mathcal{A} $ such that $p^{\prime} \in \mathcal{F} ,$ $p^{\prime} \preceq q^{\prime \prime} $ and $b $ is invertible up to pair $(p^{\prime}, q^{\prime \prime}) .$ By the same arguments as in the proof of \cite[Proposition 2.11]{KL} we may construct a net of pairs of projections $(p_{\alpha},q_{\alpha}^{\prime} )$ such that $b $ is invertible up to pair $(p_{\alpha},q_{\alpha}^{\prime}) ,$ $p_{\alpha} \geq p $ is an approximate unit for $\mathcal{F} ,$ $q_{\alpha}^{\prime}-q^{\prime} \sim p_{\alpha}-p $ and $(1-q_{\alpha})a(p_{\alpha}-p)=0 $ for all $\alpha$. Then, by the same arguments as in the proof of Lemma \ref{l 02} there exist projections $p^{\prime \prime}, \tilde{q}^{\prime \prime} $ in $\mathcal{A} $ such that $p^{\prime \prime} \sim p^{\prime} ,$ $\tilde{q}^{\prime \prime} \sim q^{\prime \prime} ,$ $b$ is invertible up to pair $(p^{\prime \prime},\tilde{q}^{\prime \prime} ) ,$ $\tilde{q}^{\prime \prime} b (1-p^{\prime \prime} )  =0$ and $p_{\alpha} \geq p^{\prime \prime} $ for all $\alpha \geq \alpha_{0}.$ Once again we can repeat the process from the proof of \cite[Proposition 2.11]{KL}  in order to construct the net of projections $\lbrace q_{\alpha}^{\prime \prime} \rbrace_{\alpha \geq \alpha_{0}} $ such that the pairs $( p_{\alpha}, q_{\alpha}^{\prime \prime} )_{\alpha \geq \alpha_{0}} $ satisfy the conditions of \cite[Proposition 2.11]{KL} with respect to the pair $( p^{\prime \prime} , \tilde{q}^{\prime \prime}) .$ In particular, $b$ is invertible up to pair $(p_{\alpha}, q_{\alpha}^{\prime \prime}) $ and $q_{\alpha}^{\prime \prime}-\tilde{q}^{\prime \prime} \sim p_{\alpha}-p^{\prime \prime} $ for all $\alpha \geq \alpha_{0} .$ Now, by Lemma \ref{r10l 1.2} we have $p^{\prime \prime} \preceq \tilde{q}^{\prime \prime}  $ since $p^{\prime \prime} \sim p^{\prime } ,$ $\tilde{q}^{\prime \prime} \sim q^{\prime \prime}  $ and $p^{\prime} \preceq q^{\prime \prime} .$ By Corollary \ref{r10c 1.5} we also get then that 
	$q_{\alpha}^{\prime \prime} = (q_{\alpha}^{\prime \prime} - \tilde{q}^{\prime \prime})+ \tilde{q}^{\prime \prime} \succeq (p_{\alpha} - p^{\prime \prime}) +p^{\prime \prime} = p_{\alpha}$ for all $\alpha \geq \alpha_{0} .$ Now, by \cite[Proposition 2.8]{KL} for each $\alpha \geq \alpha_{0} $ there exist projections $\tilde{q}_{\alpha}^{\prime }, \tilde{q}_{\alpha}^{\prime \prime } $ such that $\tilde{q}_{\alpha}^{\prime } \sim q_{\alpha}^{\prime } ,$ $\tilde{q}_{\alpha}^{\prime \prime} \sim q_{\alpha}^{\prime \prime },$ $\tilde{q}_{\alpha}^{\prime } b (1-p_{\alpha}) = \tilde{q}_{\alpha}^{\prime \prime}b (1-p_{\alpha})=0$ and $b$ is invertible both up to pair $(p_{\alpha}, \tilde{q}_{\alpha}^{\prime} )  $ and up to pair $(p_{\alpha}, \tilde{q}_{\alpha}^{\prime \prime} ) .$ By the same arguments as in the proof of Lemma \ref{l 02} we can deduce that $\tilde{q}_{\alpha} =  \tilde{q}_{\alpha}^{\prime \prime} .$ Hence $q_{\alpha}^{\prime } \sim q_{\alpha}^{\prime \prime} .$ However, we have that $p_{\alpha}  \preceq  q_{\alpha}^{\prime \prime},$ so by Lemma \ref{r10l 1.2} we get $p_{\alpha}  \preceq  q_{\alpha}^{\prime } $ for all $ \alpha \geq \alpha_{0}.$ Now we recall that $a$ is invertible up to pair $(p, q)$ and $qa(1-p)=0 .$ Since  $a(1-p)$ has left inverse up to $p$ and $p_{\alpha} \geq p $ for all $\alpha \geq \alpha_{0} ,$ there exist by \cite[Lemma 2.7]{KL}  for all $\alpha \geq \alpha_{0} $ projections $q_{\alpha} $ such that $a$ is invertile up to pair $(p_{\alpha} , q_{\alpha}) $ and $q_{\alpha} - q \sim p_{\alpha} - p .$ This follows by letting for $\alpha \geq \alpha_{0} ,$ $1-p_{\alpha} $ and $1-q_{\alpha} $ play the role of $r $ and $s $ in \cite[Lemma 2.7]{KL}, respectively. Since $q_{\alpha} - q \sim p_{\alpha} - p \sim q_{\alpha}^{\prime } - q^{\prime } $ and $q \sim q^{\prime} $ by Lemma \ref{r10l 1.2} we get that $q_{\alpha} \sim q_{\alpha}^{\prime } ,$ hence $ p_{\alpha} \preceq q_{\alpha} $ for all $\alpha \geq \alpha_{0} .$ This is a contradiction since we assumed that 
	$a \in  \mathcal{K}\Phi_{+}(\mathcal{A}) \setminus \mathcal{K}\Phi_{+}^{-} (\mathcal{A}) .$ Thus we must have that $\mathcal{K}\Phi_{+}(\mathcal{A}) \setminus \mathcal{K}\Phi_{+}^{-} (\mathcal{A}) $ is open. 
	
	The proof of the third statement is similar, whereas the second statement can be deduced from the first statement by passing to the adjoints and using Lemma \ref{l 01}.
\end{proof}

\begin{corollary}
	Let $f:[0,1] \rightarrow \mathcal{A} $ be a continuous map such that $f([0,1]) \subseteq \mathcal{K}\Phi_{+} (\mathcal{A}) \cup \mathcal{K}\Phi_{-} (\mathcal{A}) .$\\
	Then the following statements hold.\\
	1) If $f(0) \in \mathcal{K}\Phi_{+} (\mathcal{A}) \setminus \mathcal{K}\Phi (\mathcal{A}) ,$ then $f(1) \in \mathcal{K}\Phi_{+} (\mathcal{A}) \setminus \mathcal{K}\Phi (\mathcal{A}) .$\\
	2) If $f(0) \in \mathcal{K}\Phi_{-} (\mathcal{A}) \setminus \mathcal{K}\Phi (\mathcal{A}) ,$ then $f(1) \in \mathcal{K}\Phi_{-} (\mathcal{A}) \setminus \mathcal{K}\Phi (\mathcal{A}) .$\\
	3) If $f(0) \in \mathcal{K}\Phi_{+}^{-} (\mathcal{A})  ,$ then $f(1) \in \mathcal{K}\Phi_{+}^{-} (\mathcal{A})  .$\\
	4) If $f(0) \in \mathcal{K}\Phi_{-}^{+} (\mathcal{A}) ,$ then $f(1) \in \mathcal{K}\Phi_{-}^{+} (\mathcal{A})  .$\\
	5) If $f(0) \in \mathcal{K}\Phi_{0} (\mathcal{A})  ,$ then $f(1) \in \mathcal{K}\Phi_{0} (\mathcal{A})  .$\\
	6) If $f(0) \in \mathcal{K}\Phi_{+} (\mathcal{A}) \setminus \mathcal{K} \Phi_{+}^{-} (\mathcal{A})  ,$ then $ f(1) \in \mathcal{K}\Phi_{+} (\mathcal{A}) \setminus \mathcal{K}\Phi_{+}^{-} (\mathcal{A}) .$\\
	7) If $f(0) \in \mathcal{K}\Phi_{-} (\mathcal{A}) \setminus \mathcal{K} \Phi_{-}^{+} (\mathcal{A}) ,$ then $f(1) \in \mathcal{K}\Phi_{-} (\mathcal{A}) \setminus \mathcal{K} \Phi_{-}^{+} (\mathcal{A})   .$\\
	8) If $f(0) \in \mathcal{K}\Phi (\mathcal{A}) \setminus \mathcal{K} \Phi_{0} (\mathcal{A}) ,$ then $f(1) \in \mathcal{K}\Phi (\mathcal{A}) \setminus \mathcal{K} \Phi_{0} (\mathcal{A}) .$\\
	9) If $f(0) \in \mathcal{K}\Phi_{+} (\mathcal{A}) ,$ then $f(1) \in \mathcal{K}\Phi_{+} (\mathcal{A}) .$\\
	10) If $f(0) \in \mathcal{K}\Phi_{-} (\mathcal{A}) ,$ then $f(1) \in \mathcal{K}\Phi_{-} (\mathcal{A}) .$\\	
	11) If $f(0) \in \mathcal{K}\Phi (\mathcal{A}) ,$ then $f(1) \in \mathcal{K}\Phi (\mathcal{A})  $ and $index \ f(0)=index \ f(1) .$ 
\end{corollary}

\begin{proof}
	This can be shown in exactly the same way as in the proof of \cite[Corollary 4.35]{IS1}. 
\end{proof} 

\begin{corollary}
	Let $ a \in \mathcal{A}.$ Then the following statements hold.\\
	1) If $a$ belongs to the boundary of $\mathcal{K}\Phi (\mathcal{A})  $ in $\mathcal{A} ,$ then $a \in \mathcal{A} \setminus ( \mathcal{K}\Phi_{+} (\mathcal{A}) \cup \mathcal{K}\Phi_{-} (\mathcal{A}) )   .$\\
	2) If $a$ belongs to the boundary of $\mathcal{K}\Phi_{+}^{-} (\mathcal{A})  $ in $\mathcal{A} ,$ then $a \in \mathcal{A} \setminus \mathcal{K}\Phi_{+} (\mathcal{A}) .$\\
	3) If $a$ belongs to the boundary of $ \mathcal{K}\Phi_{-}^{+} (\mathcal{A}) $ in $\mathcal{A} ,$ then $a \in \mathcal{A} \setminus \mathcal{K}\Phi_{-} (\mathcal{A}) .$\\
	4) If $a$ belongs to the boundary of $\mathcal{K}\Phi_{0} (\mathcal{A})  $ in $\mathcal{A} ,$ then $ a \in \mathcal{A} \setminus \mathcal{K}\Phi (\mathcal{A}) .$
\end{corollary}

Next, we show that the set of semi-Weyl type elements is invariant under perturbations by finite type elements.  

\begin{proposition} \label{r10p 1.6}
	Let $a \in \mathcal{A}.$ Then the following holds.\\
	1) If $ a \in \mathcal{K}\Phi_{+}^{-}(\mathcal{A})$ and $f \in \mathcal{F} ,$ then $a+f \in \mathcal{K}\Phi_{+}^{-}(\mathcal{A}) .$\\
	2) If $a \in \mathcal{K}\Phi_{-}^{+}(\mathcal{A}) $ and $f \in \mathcal{F} ,$ then $a+f \in \mathcal{K}\Phi_{-}^{+}(\mathcal{A}) .$\\
	3) If $a \in \mathcal{K}\Phi_{0}(\mathcal{A}) $ and $f \in \mathcal{F} ,$ then $a+f \in \mathcal{K}\Phi_{0}(\mathcal{A}) . $
\end{proposition}

\begin{proof}
	We prove $1)$ first. Let $p, q$ be projections in $ \mathcal{A} $ such that $a$ is invertible up to pair $(p, q),$ $ p \in \mathcal{F}$ and $p \preceq q .$ In fact, by \cite[Proposition 2.8]{KL} and Lemma \ref{r10l 1.2} we may without loss of generality assume that $qa(1-p)=0 .$ By the proof of \cite[Proposition 2.11]{KL} there exists an approximate unit $\lbrace p_{\alpha} \rbrace $ for $\mathcal{F} $ and a net of projections $\lbrace q_{\alpha} \rbrace $ in $\mathcal{A}  $ such that $ p_{\alpha} \geq p ,$ $a$ is invertible up to $(p_{\alpha} , q_{\alpha}) $ and $q_{\alpha} - q \sim p_{\alpha} - p $ for all $\alpha$. \\
	Hence, by Corollary \ref{r10c 1.5} we get that $p_{\alpha} \preceq q_{\alpha} $ for all $\alpha .$ Moreover, again by \cite[Proposition 2.8]{KL} there exists a net of projections $ \lbrace \tilde{q}_{\alpha} \rbrace_{\alpha}  $ such that $a$ is invertible up to $ (p_{\alpha}, \tilde{q}_{\alpha}),$ $\tilde{q}_{\alpha} \sim q_{\alpha} $ and $ \tilde{q}_{\alpha} a (1-p_{\alpha})$ for all $ \alpha.$ Hence, by Lemma \ref{r10l 1.2} $p_{\alpha} \preceq \tilde{q}_{\alpha}  $ for all $\alpha .$ Also, we must have that $ q \leq \tilde{q}_{\alpha}$ for all $ \alpha.$ Indeed, if $b_{\alpha}$ denote $(p_{\alpha}, \tilde{q}_{\alpha}) $ inverse of $a,$ then we have 
	$$q(1-\tilde{q}_{\alpha})= q(1-\tilde{q}_{\alpha})a(1-p_{\alpha})b_{\alpha}=qa(1-p_{\alpha})b_{\alpha} = 
	qa(1-p)(1-p_{\alpha})b_{\alpha}=0, $$ 
	so $ q \leq \tilde{q}_{\alpha}$ for all $ \alpha.$ If $b$ denotes $(p, q)$ inverse of $a,$ we will prove that $ (1-p_{\alpha}) b (1-\tilde{q}_{\alpha})$ is $(p_{\alpha}, \tilde{q}_{\alpha}) $ inverse of $a.$ Indeed, 
	$$(1-p_{\alpha}) b (1-\tilde{q}_{\alpha})a(1-p_{\alpha})= (1-p_{\alpha}) b a (1-p_{\alpha})= (1-p_{\alpha}) b a (1-p) (1-p_{\alpha}) $$ 
	$$= (1-p_{\alpha}) (1-p) (1-p_{\alpha}) = 1-p_{\alpha}  .$$ 
	
	Hence, if $(1-p_{\alpha}) b_{\alpha} (1-\tilde{q}_{\alpha}	) $ is $( p_{\alpha}, \tilde{q}_{\alpha}	)$ inverse of $a$, then 
	$$(1-p_{\alpha})b_{\alpha} (1- \tilde{q}_{\alpha}) =  (1-p_{\alpha})b (1- \tilde{q}_{\alpha}) a (1-p_{\alpha}) b_{\alpha} (1-\tilde{q}_{\alpha})= (1-p_{\alpha})b (1- \tilde{q}_{\alpha})   
	, $$ 
	which proves the claim. Now, since $f \in \mathcal{F} $ and $\lbrace p_{\alpha} \rbrace $ is an approximate unit for $\mathcal{F} ,$ for $\alpha $ large enough we have 
	$\parallel f(1-p_{\alpha})\parallel \leq \parallel b \parallel^{-1} .$ Hence we get 
	$$\parallel (1-\tilde{q}_{\alpha}) f (1-p_{\alpha}) b (1-\tilde{q}_{\alpha}) \parallel \leq 1        
	\text{ and } 
	\parallel (1-p_{\alpha}) b (1-\tilde{q}_{\alpha}) f (1-p_{\alpha}) \parallel <1 $$ 
	for $ \alpha $ large enough. It follows that 
	$$1-\tilde{q}_{\alpha} + (1-\tilde{q}_{\alpha}) f (1-p_{\alpha}) b (1-\tilde{q}_{\alpha})  \text{ and } 1-p_{\alpha}+(1-p_{\alpha}) b (1-\tilde{q}_{\alpha}) f (1-p_{\alpha}) $$ 
	are invertible elements in the corner $C^{*}-$algebras 
	$$(1-\tilde{q}_{\alpha}) \mathcal{A} (1-\tilde{q}_{\alpha}) \text{ and } (1-p_{\alpha}) \mathcal{A} (1-p_{\alpha}),$$
	respectively. Hence there exists elements $c_{\alpha} $ and $\tilde{c}_{\alpha} $ in $\mathcal{A} $ such that 
	$$(  (1-\tilde{q}_{\alpha}) a (1-p_{\alpha})  +(1-\tilde{q}_{\alpha}) f (1-p_{\alpha})) (1-p_{\alpha}) b  (1-\tilde{q}_{\alpha}) \tilde{c}_{\alpha} =1-\tilde{q}_{\alpha}$$
	and
	$$c_{\alpha}  (1-p_{\alpha}) b (1-\tilde{q}_{\alpha}) (  (1-\tilde{q}_{\alpha}) a (1-p_{\alpha}) + (1-\tilde{q}_{\alpha}) f (1-p_{\alpha}) ) = 1-p_{\alpha} $$  
	where we have used that 
	$(1-p_{\alpha}) b (1-\tilde{q}_{\alpha})  $ is $(p_{\alpha} , \tilde{q}_{\alpha} ) $
	inverse of $a,$ so that 
	$$(1-\tilde{q}_{\alpha}) a (1-p_{\alpha}) b (1-\tilde{q}_{\alpha}) = 1-\tilde{q}_{\alpha} .$$ 
	By the previous arguments, we must then have that 
	$$(1-p_{\alpha}) b (1-\tilde{q}_{\alpha}) \tilde{c}_{\alpha} (1-\tilde{q}_{\alpha}) =   
	(1-p_{\alpha}) c_{\alpha} (1-p_{\alpha}) b  (1-\tilde{q}_{\alpha})$$ 
	is $ (p_{\alpha} , \tilde{q}_{\alpha} ) $ inverse of $a+f.$ Since $ p_{\alpha} \preceq \tilde{q}_{\alpha}  ,$ we must have that $a+f \in \mathcal{K}\Phi_{+}^{-}(\mathcal{A})  .$ Similarly we can prove $3).$ By passing to the adjoints and using $1)$ together with Lemma \ref{l 01} we can deduce $2).$ 
\end{proof}

Now we are able to give an algebraic characterization of semi-Weyl type elements in $\mathcal{A}$ in terms of left and right invertible elements in $\mathcal{A}.$ 

\begin{proposition} \label{r10p 1.7}
	Let $a \in \mathcal{A} .$ Then the following statements hold.\\
	1) $a \in \mathcal{K}\Phi_{+}^{-}(\mathcal{A})  $ if and only if there exist a left invertible element $b \in \mathcal{A} $ and some $ f \in \mathcal{F} $ such that $ a=b+f ,$\\
	2) $a \in \mathcal{K}\Phi_{-}^{+}(\mathcal{A}) $ if and only if there exist a right invertible element $ b \in \mathcal{A} $ and some $ f \in \mathcal{F} $ such that $a=b+f ,$\\
	3) $a \in \mathcal{K}\Phi_{0} (\mathcal{A}) $ if and only if there exist an invertible element $b \in \mathcal{A}  $ and some $ f \in \mathcal{F} $ such that $ a=b+f .$ 
\end{proposition} 

\begin{proof}
	We prove 1) first. Suppose that $a$ is invertible up to $(p, q)$ where $ p \in \mathcal{F} $ and $ p\preceq q .$ By \cite[Proposition 2.8]{KL} and Lemma \ref{r10l 1.2} we may assume that $q a (1-p) = 0  .$ Let $c$ be $(p, q)$ inverse of $a$ and $v \in \mathcal{A} $ be such that $v^{*}v=p $ and $ vv^{*}=q^{\prime}$ for some projection $q^{\prime} \leq q .$ 
	Set $ b= a (1-p) + q^{\prime}vp .$ Then we get that $ (c (1-q) + v^{*}q^{\prime} ) b = c a (1-p) + v^{*} q^{\prime} v p = (1-p) + v^{*} v v^{*} v p = 1-p + p=1 $ (where we have used that $qa(1-p)=0 $ and $q^{\prime} \leq q ) .$ Hence $b$ is left invertible in $ \mathcal{A} .$ Moreover $ (a- q^{\prime}v)p \in \mathcal{F} $ since $ p \in \mathcal{F} $ and $ \mathcal{F}$ is an ideal, and $a=b+ (a - q^{\prime}v)p .$ 
	
	Conversely, if $b$ is left invertible in $ \mathcal{A},$ by \cite[Lemma 2.3]{KL} there exists a projection $s \in \mathcal{A} $ such that $sb=b $ and $b$ is invertible up to $(0, s).$ Hence $b \in \mathcal{K}\Phi_{+}^{-}(\mathcal{A}) ,$ so by Proposition \ref{r10p 1.6} we must have that $b+f \in \mathcal{K}\Phi_{+}^{-}(\mathcal{A}) $ for all $ f \in \mathcal{F} .$ 
	
	The proof of the third statement is similar, whereas the second statement can be deduced from the first statement by passing to the adjoints and using Lemma \ref{l 01}.
\end{proof}

\begin{corollary}
	The sets $\mathcal{K}\Phi_{+}^{-} (\mathcal{A}) , \mathcal{K}\Phi_{-}^{+} (\mathcal{A})  $ and $\mathcal{K}\Phi_{0} (\mathcal{A})  $ are semigroups under the multiplication. 
\end{corollary}

The next lemma will be used later in this paper in connection with von-Neumann algebras in order to show that if $\mathcal{A}$ is a properly infinite von-Neumann algebra, then the set of those Fredholm type elements that are both upper semi-Weyl and lower semi-Weyl is exactly the set of Weyl type elements in $\mathcal{A}$. 

\begin{lemma}
	Let $a \in \mathcal{K}\Phi_{+}^{-} (\mathcal{A}) \cap \mathcal{K}\Phi_{-}^{+} (\mathcal{A}) \cap \mathcal{K}\Phi (\mathcal{A}) .$ Then there exist projections $p, q$ in $\mathcal{F}$ such that $a$ is invertible up to $(p, q), qa (1-p)=0,$ $p \preceq q $ and $q \preceq p.$
\end{lemma}

\begin{proof}
	Let $p,q,p^{\prime},q^{\prime}  $ be projections in $\mathcal{A} $ such that $p, q^{\prime} \in \mathcal{F},p \preceq q , q^{\prime} \preceq p^{\prime}  $ and $a$ is invertible both up to $(p, q)$ and $(p^{\prime},q^{\prime}) .$ Since $a \in  \mathcal{K}\Phi (\mathcal{A}) ,$ by Lemma \ref{l 02} we must have that $q, p^{\prime} \in \mathcal{F} .$ By the proof of Lemma \ref{l 02} there exists an approximate unit $\lbrace p_{\alpha} \rbrace $ and a net of projections $\lbrace \tilde{q}_{\alpha} \rbrace $ in $\mathcal{A}$ and a projection $p^{\prime \prime} \sim p^{\prime}$ such that $p \leq p_{\alpha}, p^{\prime \prime} \leq p_{\alpha} ,$ $a$ is 
	invertible up to $(p_{\alpha},  \tilde{q}_{\alpha}) $ and $ \tilde{q}_{\alpha} a (1-p_{\alpha})= 0 $ for all $\alpha .$ Moreover, for each $ \alpha $ there exists projections $q_{\alpha} $ and $q_{\alpha}^{\prime \prime} $ such that $\tilde{q}_{\alpha} \sim q_{\alpha} \sim  q_{\alpha}^{\prime \prime},$ $q \leq q_{\alpha},$ $ q^{\prime} \sim q^{\prime \prime} \leq q_{\alpha}^{\prime \prime}  $ for some projection $q^{\prime \prime} $ in $\mathcal{A} $ where $p_{\alpha} - p \sim q_{\alpha} - q $ and $p_{\alpha} - p^{\prime \prime} \sim q_{\alpha}^{\prime \prime} - q^{\prime \prime} .$ Since $p \preceq q $ and $  q^{\prime \prime} \sim q^{\prime} \preceq p^{\prime} \sim p^{\prime \prime} ,$ by Corollary \ref{r10c 1.5} we get that $p_{\alpha} \preceq q_{\alpha} \sim \tilde{q}_{\alpha} $ and $p_{\alpha} \succeq q_{\alpha}^{\prime \prime} \sim \tilde{q}_{\alpha} $ for all $\alpha .$
\end{proof}

In \cite{DDj2} \DJ{}or\dj{}evi\'{c} introduced the notion of generalized Weyl operator on a Hilbert space as a closed range bounded operator whose kernel is isomorphic to the orthogonal complement of the image. In \cite[Theorem 1]{DDj2},  he proved that if a composition of two such operators has closed image, then this composition is also a generalized Weyl operator. The next proposition provides a generalization of this result from the special case when $\mathcal{A}=B(H)$ to the case of general unital $C^{*}$-algebras.

\begin{proposition} \label{R8 p2.23}
	Let, $a,b \in \mathcal{A} $ and suppose that there exist projections $p,q,p^{\prime}, q^{\prime}, \tilde{p}, \tilde{q} $ in $ \mathcal{A}$ such that $(1-q)a(1-p)=a ,$ $(1-q^{\prime})b(1-p^{\prime})=b  ,$  $(1-\tilde{q})ba(1-\tilde{p})=ba $ and $a,b,ba $ are invertible up to $(p,q), $  $(p^{\prime},q^{\prime}) $ and $(\tilde{p},\tilde{q}),$ respectively. If $p \sim q $ and $p^{\prime} \sim q^{\prime},$ then  $\tilde{p} \sim\tilde{q} .$ Moreover, $p \leq \tilde{p} ,$ $q^{\prime} \leq \tilde{q},$ $\tilde{p}-p  \preceq p^{\prime}$ and $\tilde{q} - q^{\prime}  \preceq q  .$
\end{proposition}

\begin{proof}
	Let $\tilde{c} $ be 
	$( \tilde{p}, \tilde{q}) $ inverse of $ba .$ Then 
	$$(1-\tilde{p} )p = \tilde{c}  (1-\tilde{q}) ba (1-\tilde{p}) p = \tilde{c}bap = \tilde{c} b (1-q) a (1-p ) p =0 .$$

	Similarly, we have 
	$$q^{\prime} (1-\tilde{q}) = q^{\prime}  (1-\tilde{q}) ba (1-\tilde{p}) \tilde{c} = q^{\prime} ba (1-\tilde{p}) \tilde{c} = q^{\prime} (1-q^{\prime}) b a (1-\tilde{p}) \tilde{c} =0  ,$$ 
	so $ p \leq \tilde{p} $ and $q^{\prime} \leq \tilde{q} .$ Hence, by  \cite[Lemma 2.5]{KL},  there exists a projection $t \in \mathcal{A} $ such that $t a ( \tilde{p}-p) = a ( \tilde{p}-p), $ $t \sim \tilde{p}-p $ and $a$ is invertible up to $(\tilde{p}-p,t) .$ Let, $  \tilde{a}  $ be $ (\tilde{p}-p,t)$ inverse of $a.$ Then we get that 
	$$bt=bta(\tilde{p}-p)\tilde{a}= ba (\tilde{p}-p) \tilde{a} = (1-\tilde{q})  ba (1 - \tilde{p})  ( \tilde{p} - p) \tilde{a} = 0.$$ 
	However, if $b^{\prime} $ is $(p^{\prime},q^{\prime}) $  inverse of $b,$ then we have that $$ (1-p^{\prime})t = b^{\prime}(1-q^{\prime}) b (1-p^{\prime}) t = b^{\prime} b t =0,$$
	so $t \leq p^{\prime} .$ Thus, $\tilde{p} - p \sim t \leq p^{\prime} ,$ so $ \tilde{p}-p \preceq p^{\prime} .$ Similarly we can show that $\tilde{q} - q^{\prime} \preceq q .$
	
	Now, since $b a $ is invertible up to $(\tilde{p}, \tilde{q}) $ by assumption, it follows that $a$ is left invertible up to $ \tilde{p}$ and $b$ is right invertible up to $\tilde{q} .$ By  \cite[Lemma 2.5]{KL} there exists a projection $r \in \mathcal{A} $ such that $(1-r) a (1- \tilde{p}) = 1- \tilde{p} $ and such that $a$ is invertible up to $(\tilde{p} , r) .$ Since $1-\tilde{p} \leq 1-p, $ by \cite[Lemma 2.7]{KL} we have that $1-r \leq 1-q $ (which gives $q \leq r $) and $r-q \sim \tilde{p} -p .$ This is because $qa(1-p) = q (1-q) a (1-p) = 0, $ so that we can indeed apply \cite[Lemma 2.7]{KL}. Since $ p \sim q,$ by Lemma  \ref{r10l 1.2} we get that $\tilde{q} \sim r .$ Next, since $b$ is right invertible up to $\tilde{q},$ it follows that $b^{*}$ is left invertible up to $\tilde{q}.$  Moreover, by Lemma \ref{l 01} we have that $b^{*}$ is invertible up to $(q^{\prime}, p^{\prime}) $ an $p^{\prime} b^{*}(1-q^{\prime})=0 $ since $b$ is invertible up to $(p^{\prime},q^{\prime}) $ 	and $(1-q^{\prime}) b (1-p^{\prime}) = b .$ Hence we can apply the preceding arguments on $b^{*}$ to deduce that there exists some projection $s \in \mathcal{A} $ such that $(1-s) b^{\prime} (1-\tilde{q}) = b^{\prime} (1-\tilde{q}) ,$ $p^{\prime} \leq s, b^{*}$ is invertible up to $(\tilde{q}, s) $ and $s- p^{\prime} \sim \tilde{q} - q^{\prime}  $ which gives that $ s \sim \tilde{q}.$ It follows that $b$ is invertible up to $(s,  \tilde{q} ) .$\\
	We will prove now that $b$ is also invertible up to $(r, \tilde{q} ) .$ Indeed, if $c$ is $(\tilde{p} , r) $ inverse of $a$ and $d$ is $(\tilde{p}, \tilde{q}) $ inverse of $ba,$ then we get that

	$$ a (1-\tilde{p}) d (1-\tilde{q}) b (1-r) = a (1-\tilde{p}) d (1-\tilde{q}) b (1-r) a (1-\tilde{p}) c     $$
	$$= a (1-\tilde{p}) d (1-\tilde{q}) b a (1-\tilde{p}) c = a (1-\tilde{p}) c =1-r, $$
	and 
	$$(1-\tilde{q}) b (1-r) a (1-\tilde{p}) d = (1-\tilde{q}) b a (1-\tilde{p}) d = 1-\tilde{q} ,$$ 
	so $ a (1-\tilde{p}) d  $ is $(r, \tilde{q}) $ inverse of $b.$ Since $b$ is invertible both up to $(r, \tilde{q}) $ and up to $(s, \tilde{q}) ,$ by passing to the adjoints and applying Lemma \ref{l 01} together with Corollary \ref{cor 3.3} we deduce that $ r \sim s.$ However, we recall from previous calculations that $\tilde{p} \sim r $ and $s \sim \tilde{q} .$ Hence we must have that $\tilde{p} \sim \tilde{q} ,$ and the proof is complete.
	
\end{proof} 

\section{Applications to von Neumann algebras} 

We start with the following corollary.

\begin{corollary}
	Let $\mathcal{A} $ be a von Neumann algebra acting on a Hilbert space $H$ and $T,S \in \mathcal{A} $ such that $Im T, Im S $ and $Im ST $ are closed. Then there exist closed subspaces $H_{1} $ and $H_{2} $ of $H$ such that $P_{H_{1}}, $ $P_{H_{2}}$ belong to $\mathcal{A}$,
	$$\ker ST = \ker T \oplus H_{1}, \text{ } Im ST^{\perp} = Im S^{\perp} \oplus  H_{2},  \text{ } P_{H_{1}} \preceq P_{\ker S} \text{ and } P_{H_{2}} \preceq P_{Im T^{\perp}}$$
	(where $P$ stands for orthogonal projection onto the respective subspace). Moreover, if $ P_{\ker T} \sim P_{Im T^{\perp}} $ and $P_{\ker S} \sim P_{Im S^{\perp}} ,$ then $P_{\ker ST} \sim P_{Im ST^{\perp}} .$
\end{corollary}

\begin{proof}
	The statements follow from Proposition \ref{R8 p2.23}. 
\end{proof}

In the sequel $\mathcal{A}$ denotes a properly infinite von Neumann algebra and $Proj_{0}(\mathcal{A})$ denotes the set of all finite projections in $\mathcal{A},$ (i.e. those projections that are not Murray von Neumann equivalent to any of its subprojections). \\
We recall the notion of $\mathcal{A}$-Fredholm operator, originally introduced by Breuer in \cite{BR}, \cite{BR2}. 

\begin{definition} \cite[Definition 3.1]{KL}
	The operator $T \in \mathcal{A}$ is said to be $\mathcal{A}$-Fredholm if the following holds.\\
	$(i)$  $P_{\ker T}  \in Proj_{0}(\mathcal{A}), $ where $P_{\ker T}$ is the projection onto the subspace $\ker T$.\\
	$(ii)$ There is a projection $E \in Proj_{0}(\mathcal{A}) $ such that $Im(I - E) \subseteq Im T . $\\
	The second condition ensures that $P_{\ker T^{*}}$
	also belongs to $Proj_{0}(\mathcal{A}).$\\
	
	The index of an $\mathcal{A}$-Fredholm operator $T$ is defined as 
	$$index T = dim(\ker T ) - dim(\ker T^{*}) \in I(\mathcal{A}) . $$
	
	Here, $I(\mathcal{A})$ is the so called\textit{ index group} of a von Neumann algebra $\mathcal{A}$ defined as the Grothendieck group of the commutative monoid of all representations of the commutant $\mathcal{A}^{\prime}$ generated by representations of the form $\mathcal{A}^{\prime}  \ni S \mapsto ES =\pi_{E}(S) $ for
	some $E \in Proj_{0}(\mathcal{A}) .$ For a subspace $L,$ its dimension $dim L$ is defined as the $class [\pi_{P_{L}}] \in I(\mathcal{A}) $ of the representation $ \pi_{ P_{L}},$ where $P_{L}$ is the projection onto $L.$
\end{definition}

We give then the following characterization of $\mathcal{A}$-Fredholm operators. 

\begin{lemma}  \label{r12 l11}
	Let $\mathcal{A} $ be a properly infinite von Neumann algebra. Then an operator $T \in \mathcal{A} $ is $\mathcal{A}-$Fredholm in the sense of Breuer if and only if there exist projections $ P,Q \in Proj_{0} (\mathcal{A})$ such that $T$ is invertible up to $(P,Q) .$ 
\end{lemma}

\begin{proof}
	By the proof of \cite[Corollary 3.10]{KL}  we may let $\mathcal{F} = \mathfrak{m} $ where $\mathfrak{m}$ is the norm closure of the set of all $S \in \mathcal{A} $ for which $P_{\overline{ImS}} \in Proj_{0} (\mathcal{A}) .$ Moreover, in that proof it has also been shown that $Proj_{0} (\mathcal{A})$ is an approximate unit for $ \mathfrak{m}.$ Hence, by \cite[Lemma 2.9]{KL} if $P \in \mathfrak{m} $ and $P$ is a projection, then $ P \preceq P_{\alpha_{0}} $ for some $P_{\alpha_{0}} \in Proj_{0} (\mathcal{A}).$ It follows then that $P \in Proj_{0} (\mathcal{A}).$ Hence 
	\begin{equation} \label{f1}
		\lbrace P \in \mathfrak{m} \text{ } \vert \text{ } P \text{   is projection  } \rbrace = Proj_{0}(\mathcal{A}). 
	\end{equation}
	Moreover, by \cite[Corollary 3.10]{KL} an operator $T \in \mathcal{A} $ is $\mathcal{A}-$Fredholm in the sense of Breuer if and only if $T$ is invertible up to $(P,Q) $ for some 
	projections $P,Q \in \mathfrak{m} .$ Combining all these facts together, we obtain the lemma.
\end{proof}

The relation (\ref{f1}) has several useful applications in semi-$\mathcal{A}$-Fredholm theory, as listed below.

\begin{definition} \label{r12 d15}
	Let $\mathcal{A} $ be a properly infinite von Neumann algebra and $T \in \mathcal{A} .$ We say that $T$ is upper semi$-\mathcal{A}-$Fredholm if there exist projections $P, Q$ in $\mathcal{A}$ such that $T$ is invertible up to $(P,Q)$ where $ P \in Proj_{0} (\mathcal{A}).$ If in addition $P \preceq Q ,$ we say that $T$ upper semi$-\mathcal{A}-$Weyl. Similarly we say that $T$ is lower semi$- \mathcal{A} -$Fredholm an lower semi$-\mathcal{A} -$Weyl, however in this case we assume that $Q \in  Proj_{0} (\mathcal{A})$ and $Q \preceq P .$
\end{definition}

\begin{corollary} \label{r12 c17}
	Let $T \in \mathcal{A}.$ Then $T$ is upper (respectively lower) semi-Fredholm type element in $ \mathcal{A}$ with respect to $\mathfrak{m}$  if and only if $T$ is upper (respectively lower) semi$- \mathcal{A} -$Fredholm. Similarly, $T$ is upper (respectively lower) semi-Weyl type element in $\mathcal{A} $ with respect to $\mathfrak{m}$  if and only if $T$ is upper (respectively lower) semi$- \mathcal{A} -$Weyl. Finally, $T$ is Weyl type element in $ \mathcal{A} $ with respect to $\mathfrak{m}$  if and only if $T$ is $\mathcal{A} -$Weyl.
\end{corollary}

\begin{proof}
	By (\ref{f1}) we have that $ \mathfrak{m} \cap Proj(\mathcal{A})= Proj_{0}(\mathcal{A}) .$ 
\end{proof}

\begin{corollary}
	Let $T \in \mathcal{A} .$ Then $T$ is upper semi-$\mathcal{A} $-Fredholm if and only if there exists some $P \in Proj_{0} (\mathcal{A}) $ such that $T$ is bounded below on $(I-P)(H) .$ Similarly, $T$ is lower semi- $\mathcal{A}$-Fredholm if and only if there exists some $ Q \in Proj_{0} (\mathcal{A}) $ such that $(I-Q)(H) \subseteq Im T .$
	
\end{corollary}

\begin{proof}
	Combine Lemma \ref{l 06} together with Corollary \ref{r12 c17} and  (\ref{f1}).
\end{proof}

\begin{corollary}
	Let $T \in \mathcal{A} .$ Then the following statements hold.\\
	1) If $T$ is upper semi$- \mathcal{A} -$Fredhoolm, then $ P_{\ker T}  \in Proj_{0} ( \mathcal{A}).$ In particular, if $Im T$ is closed, then $T$ is upper semi$- \mathcal{A} -$Fredholm if and only if $P_{\ker T}  \in Proj_{0} ( \mathcal{A}) .$\\
	2) If $T$ is lower semi$- \mathcal{A}-$Fredholm, then $P_{\overline{Im T}} \in  Proj_{0} ( \mathcal{A}) .$ In particular, if $Im T$ is closed, then $T$ is lower semi$- \mathcal{A} -$Fredholm if and only if $ P_{Im T}  \in Proj_{0} ( \mathcal{A}).$
\end{corollary}

\begin{proof}
	Apply Lemma \ref{l 06} together with (\ref{f1}). 
\end{proof}

\begin{corollary}
	We have that $\mathcal{K}\Phi_{+}^{-} (\mathcal{A}) \cap \mathcal{K}\Phi_{-}^{+} (\mathcal{A}) \cap \mathcal{K}\Phi (\mathcal{A}) = \mathcal{K}\Phi_{0} (\mathcal{A}) .$ 
\end{corollary}

\begin{proof}
	By (\ref{f1}) we have that if $P, Q$ are projections in $ \mathfrak{m},$ then $P,Q \in Proj_{0}(\mathcal{A}) .$ Now, if $P \preceq Q $ and $ Q \preceq P$ for some	$P,Q \in Proj_{0}(\mathcal{A}) ,$ by \cite{BR}, \cite{BR2} it follows that $P \sim Q .$ 
\end{proof}

\begin{corollary}
	Let $T \in \mathcal{A}.$ Then the following statements hold.\\
	1) $T$ is upper semi-$\mathcal{A}$- Weyl if and only if there exist some $S \in \mathcal{A} $ and $F \in \mathfrak{m} $ such that $S$ is bounded below and $T=S+F.$\\
	2) $T$ is lower semi-$\mathcal{A}$- Weyl if and only if there exist some $S \in \mathcal{A} $ and $F \in \mathfrak{m} $ such that $S$ is surjective and $T=S+F.$\\
	3) $T$ is $\mathcal{A}$-Weyl if and only if there exist some $S \in \mathcal{A}$ and $F \in\mathfrak{m} $ such that $S$ is invertible and $T=S+F.$
\end{corollary}

\begin{proof}
	The statements follow from Proposition \ref{r10p 1.7}.
\end{proof}

\begin{lemma}  \label{r102 l3.6} 
	Let $T \in \mathcal{A} .$ Suppose that $T$ is $\mathcal{A}$-Fredholm and that  $Im \ T$ is closed. Then there exist an $\epsilon > 0 $ such that for every $ S \in \mathcal{A}$ with $\parallel S \parallel < \epsilon ,$ we have that $P_{\ker (T+S)} \preceq P_{\ker T} $ and $P_{Im(T+S)^{\perp}} \preceq P_{Im T^{\perp}}.$
\end{lemma} 

\begin{proof}
	Since $T$ is $\mathcal{A} $-Fredholm, by Lemma \ref{l 06},  Corollary \ref{r12 c17} and  (\ref{f1}) we must have that $P_{\ker T},P_{ImT^{\perp}} \in Proj_{0} (\mathcal{A}) .$ Moreover, $T$ is invertible up to $(P_{\ker T},P_{ImT^{\perp}})$ since $Im T$ is closed. By \cite[Lemma 2.4]{KL} there exists an $\epsilon > 0  $ such that if $S \in  \mathcal{A}$ and $\parallel S \parallel < \epsilon ,$ then $T+S$ is invertible up to $(P_{\ker T},P_{ImT^{\perp}}) .$ Hence, by \cite[Proposition 2.8]{KL} there exist some $\tilde{P}, \tilde{Q} \in Proj_{0}(\mathcal{A}) $ such that $\tilde{P} \sim P_{\ker T}, \tilde{Q} \sim P_{ImT^{\perp}}, T+S $ is invertible both up to $(\tilde{P} , P_{ImT^{\perp}}) $ and $(P_{\ker T}, \tilde{Q}) ,$ and $\tilde{Q}(T+S)(I-P_{\ker T}) = (I-P_{Im T^{\perp}})(T+S) \tilde{P} = 0 .$ By the proof of Lemma \ref{l 06} we must have that $P_{\ker (T+S)}  \leq \tilde{P} $ and $P_{Im (T+S)^{\perp}}   \leq \tilde{Q} $ which proves the lemma.
\end{proof}

The next proposition is a generalization of the punctured neighborhood theorem \cite[ Theorem 1.7.7]{ZZRD} in the setting of von Neumann algebras. The main idea in this proof is motivated by the proof of \cite[Theorem 3.26]{IS3}. 

\begin{proposition}
	Let $\mathcal{A} $ be a properly infinite von Neumann algebra acting on a Hilbert space $H$ and $T$ be an $\mathcal{A} $-Fredholm operator. Suppose that $Im T^{n} $ is closed for all $n$ and that 
	$$T (\cap_{n=1}^{\infty} \ Im T^{n}) = \cap_{n=1}^{\infty} Im T^{n} .$$ 
	Then there exists an $ \epsilon >0$ such that if $\lambda \in \mathbb{C} $ and $0 < \vert \lambda \vert < \epsilon,$ then $dim \ \ker T = dim \ \ker(T-\lambda I)+ dim \ N_{1}  $ and $dim \ (Im T)^{\perp} = dim \ (Im(T-\lambda I))^{\perp} + dim \ N_{1} $ in $I(\mathcal{A}) $ for some fixed closed subspace $N_{1}$ of $H$ where $P_{N_{1}} \in Proj_{0} (\mathcal{A}). $
\end{proposition}

\begin{proof}
	Since $T$ is $\mathcal{A} $-Fredholm and $Im T$ is closed, by Lemma \ref{r102 l3.6}, \cite[Corollary 3.10]{KL} and \cite[Proposition 2.13]{KL} there exists an $\epsilon_{1} > 0  $ such that if $\lambda \in \mathbb{C} $ and $\vert \lambda \vert < \epsilon_{1} ,$ then $P_{\ker (T-\lambda I)} \preceq P_{\ker T}, $ $P_{Im (T-\lambda I)^{\perp}} \preceq P_{(Im T)^{\perp}} $ and $index \ (T-\lambda I)= index \ T .$ Also, we have that 
	$$\ker (T-\lambda I) \subseteq Im^{\infty} (T) := \cap_{n=1}^{\infty} Im(T^{n}) $$
	when $\lambda \neq 0 .$ 
	Since 
	$$P_{Im^{\infty}(T)}= s - \lim_{n \rightarrow \infty} P_{Im (T^{n})} ,$$ 
	we have $P_{Im^{\infty}(T)} \in \mathcal{A} .$ Furthermore, 
	$$( \ker T \cap Im^{\infty}(T) ) \oplus (\ker T)^{\perp} =  \ker  P_{Im^{\infty}(T)^{\perp}} P_{\ker T} ,$$ 
	so $P_{\ker T \cap Im^{\infty}(T) \oplus (\ker T)^{\perp}}$	belongs to $ \mathcal{A}.$ Hence $P_{\ker T \cap Im^{\infty}(T) } $ which is equal to the diference $$P_{\ker T \cap Im^{\infty}(T) \oplus (\ker T)^{\perp}} - P_{(\ker T)^{\perp}} ,$$ 
	also belongs to $\mathcal{A} .$ Moreover, $ P_{\ker T \cap  Im^{\infty}(T) } \leq P_{\ker T},$
	so we must have that $P_{\ker T \cap Im ^{\infty} (T)}  \in Proj_{0} (\mathcal{A})$ as $P_{\ker T } \in  Proj_{0} (\mathcal{A}) .$ Set $M=\bigcap_{n=1}^{\infty} Im (T^{n}) .$ Then $T(M)=M$ by the assumption. Consider the operator $ TP_{M}+I-P_{M}.$ Then $Im (TP_{M}+I-P_{M})=H $ and $\ker (TP_{M} + I-P_{M}) = \ker T \cap M .$ Since $P_{\ker T \cap M}  \in Proj_{0} (\mathcal{A}) ,$ we get that $TP_{M}+I-P_{M} $ is an $\mathcal{A}$-Fredholm operator. Put $S:= TP_{M}+I-P_{M} .$ Since $S$ is surjective, by the same arguments as above there exists an $ \epsilon_{2} >0$ such that whenever $\lambda \in \mathbb{C} $ and $\vert \lambda \vert < \epsilon_{2} ,$ we have that $P_{\ker (S-\lambda I)}  \preceq P_{\ker S},$ $P_{Im (S-\lambda I)^{\perp}}  \preceq  P_{Im S^{\perp}   } $ and $index \ (S-\lambda I) = index S.$ Since $Im S^{\perp} = \lbrace 0 \rbrace ,$ we must hence have that $Im (S-\lambda I)^{\perp} = \lbrace 0 \rbrace $ for all $\lambda \in \mathbb{C} $ with $ \vert \lambda \vert < \epsilon_{2} .$ Since $S- \lambda I $ has the matrix 
	$ 
	\begin{pmatrix}
		T - \lambda 1 & 0 \\
		0 & 1 - \lambda  
	\end{pmatrix} 
	$  
	with respect to the decomposition $M \oplus M^{\perp} $ and $\ker (T-\lambda I) \subseteq M ,$ we get that $\ker (S-\lambda I) = \ker (T-\lambda I)$ for $0 < \vert \lambda \vert < 1 .$ Therefore, $dim \ \ker (T-\lambda I) = dim \ \ker (S- \lambda I) =index(S-\lambda I)=index  S = dim \ \ker S$ whenever $ 0 < \vert \lambda \vert < min \lbrace \epsilon_{2} , 1 \rbrace.$ Moreover, since $S$ has the matrix 
	$ 
	\begin{pmatrix}
		T  & 0 \\
		0 & 1  
	\end{pmatrix} 
	$   
	with respect to the decomposition $M \oplus M^{\perp} ,$ we get that $ \ker S = \ker T \cap M.$ \\
	Let $N_{1} $ be the orthogonal complement of $ \ker T \cap M$ in $\ker T.$ Then $P_{N_{1}} = P_{\ker T} - P_{\ker T \cap M} , $ so $ P_{N_{1}} \in \mathcal{A}.$ Furthermore, $P_{N_{1}} \leq P_{\ker T} ,$ so $P_{N_{1}} \in Proj_{0}(\mathcal{A}) $ as $P_{\ker T} \in Proj_{0} (\mathcal{A}) .$ Therefore, if $0 < \vert \lambda \vert < min \lbrace \epsilon_{2} , 1 \rbrace ,$ then 
	$$dim \ \ker T = dim \ \ker T \cap M + dim \ ker N_{1} = dim  \ker S + dim N_{1}$$ 
	$$= dim \ \ker (S-\lambda I) + dim N_{1}  
	= dim \ker (T-\lambda I) + dim N_{1} $$ 
	whenever $ 0 < \vert \lambda \vert < min \lbrace \epsilon_{2} , 1 \rbrace .$ If in addition $ \vert \lambda \vert < \epsilon_{1} ,$ then $index \ (T-\lambda I) = index \ T .$ So, if $  0 < \vert \lambda \vert < min \lbrace \epsilon_{1}, \epsilon_{2} ,  1 \rbrace,$ then $index \ (T-\lambda I) = index \ T $ and $dim \ker T = dim \ker (T - \lambda I) + dim N_{1} .$ Hence 	
	$$ \dim (Im T)^{\perp} = dim  (Im (T-\lambda I))^{\perp} + dim N_{1}.$$
\end{proof}


\bibliographystyle{amsplain}

\vspace{.1in}
\end{document}